\documentclass[a4j,10pt]{article}
\usepackage{amsmath}
\usepackage{amssymb}
\usepackage{amsthm}
\usepackage{enumerate}
\newtheorem{thm}{Theorem}[section]
\newtheorem{prop}{Proposition}[section]
\newtheorem{lem}{Lemma}[section]
\newtheorem{cor}{Corollary}[section]
\newtheorem{ex}{Example}[section]
\newtheorem{rem}{Remark}[section]

\begin{document}
\title{An analogue of Girstmair's formula in function fields }
\author{By D. Shiomi}
\date{}
\maketitle

%
%
%
%
\section{Introduction}
Suppose that  $p$ is an odd prime 
and $g>1$ is  a primitive root modulo $p$. 
Let $M$ be a number field 
contained in the $p$-th cyclotomic field. 
In \cite{gi1, gi2}, Girstmair found  a surprising relation 
between 
the relative class number of $M$
and the digits of $1/p$ in base $g$. 
In particular, he obtained 
%
%
\begin{thm}
\label{th11}
Suppose that $p\ge 7$ and $p \equiv 3 \pmod 4$. Then, 
%
%
\begin{eqnarray*}
(g+1)h(-p)=\sum_{k=1}^{p-1}(-1)^kx_k, 
\end{eqnarray*}
 where $h(-p)$ is the class number of 
$\mathbb{Q}(\sqrt{-p})$, and
$
1/p=\sum_{k=1}^{\infty}x_kg^{-k}
$
is the digit expansion of $1/p$ in base $g$. 
\end{thm}
After his works, this type of formulas has been 
further studied  by several authors (cf.  \cite{hi, kk, mi,mt}). 
%
%
%

%
%
In this paper, we consider an analogue of
Girstmair's formula in function fields. 
Let $\mathbb{F}_q$ be the finite field with $q$ elements of characteristic $p$. 
Let $K=\mathbb{F}_q(T)$ be the rational function field over $\mathbb{F}_q$,
and  put $\mathbb{A}=\mathbb{F}_q[T]$. 
Suppose that 
 $P \in \mathbb{A}$ is monic irreducible 
 and $G \in \mathbb{A}$ is a primitive root modulo $P$.  
Let $L$ be an  extension field of $K$
contained in the $P$-th cyclotomic function field. 
We denote by $h_L$  the divisor class number of $L$, i.e., 
the order of the divisor class group of $L$ of degree $0$. 
 
The goal of this paper is to give relations between the
plus and minus parts of $h_L$
and the digits of $1/P$ in base $G$
(see Theorem \ref{th31} and \ref{th32}).  
As a consequence, we prove 
the following theorem,
which can be regarded as an analogue of Theorem \ref{th11}. 
%
%
%
%
\begin{thm}
\label{th12}
Suppose that $\gcd (2,q)=1$, $d=\deg P$, and 
$\deg G \ge \deg P$. 
\vspace{1mm}
Let $L=K(\sqrt{(-1)^dP})$. Let
$
1/P=\sum^{\infty}_{k=1}H_k/G^k
$
be the digit expansion of $1/P$ in base $G$. 
Then, 
\begin{eqnarray*}
 h_L=
 \begin{cases}
\;  \displaystyle\sum_{k=1}^{r}(-1)^k\deg H_k
 \hspace{15mm} \text{if $d$ is even, } \\[5mm]
(-\eta_G) \displaystyle\sum_{k=1}^{r}(-1)^k\varepsilon_k, 
 \hspace{12mm}\text{if $d$ is odd}.  
 \end{cases}
\end{eqnarray*}
Here, $r=(q^d-1)/(q-1)$ and  
\begin{eqnarray*}
\eta_G=
 \begin{cases}
 1
 \hspace{6mm} 
 \text{if $\delta(G)\in (\mathbb{F}_q^{\times})^2$,  }  \\[3mm]
 -1
  \hspace{3mm}\text{if $\delta(G)\not\in (\mathbb{F}_q^{\times})^2$, }
  \end{cases}
  \hspace{2mm}
\varepsilon_k=
 \begin{cases}
 1
 \hspace{8mm} 
 \text{if $\delta(H_k)\in (\mathbb{F}_q^{\times})^2$,  }  \\[3mm]
 -1
  \hspace{5mm}\text{if $\delta(H_k)\not\in (\mathbb{F}_q^{\times})^2$}, 
  \end{cases}  
\end{eqnarray*}
where $\delta(I)$ 
is the leading coefficient of $I \in \mathbb{A}$ if $ I\neq 0$, 
and $\delta(0)=0$. 
 \end{thm}
%
%
\begin{rem}
Assume the notation in Theorem \ref{th12}. 
Consider 
\begin{eqnarray*}
S_{P,G}=\sum_{k=1}^{q^d-1}(-1)^kH_k,
\end{eqnarray*}
which is analogous to the alternating sum in Theorem \ref{th11}. 
However, this sum seems to have no connection 
with the divisor class number of $K(\sqrt{(-1)^dP})$.
In fact, we have $S_{P,G}=0$ by 
considering the case $\alpha=-1$ in 
Corollary \ref{cor21} (see Section 2.1). 
\end{rem}

This paper is organized as follows. In Section 2.1, 
we study about
some properties of 
digits of  rational functions over finite fields.
In Section 2.2, we review  basic facts about
cyclotomic function fields and some formulas for divisor class numbers. 
In Section 3,  we give the main results of this paper
and some examples. 
%
%
%
%
\section{Preparations}
\subsection{The digit expansion}
Let $v_\infty$ be the valuation of $K$ satisfying with $v_{\infty}(T) =-1$. 
Then 
\[
v_{\infty}(F_1/F_2)=\deg F_2-\deg F_1 \hspace{5mm}  
(F_1, F_2 \in \mathbb{A} \; \backslash \; \{0\}).
\]
 We denote by $K_\infty$ the completion of $K$ by $v_{\infty}$. 
%
%
%
%
\begin{prop}
\label{prop21}
Let  $G \in \mathbb{A}$ with $\deg G \ge 1$.
Define
\[
S_G=\{I \in \mathbb{A}\; | \deg I<\deg G\} \cup \{0\}.
\]
Any rational function $f \in K$
can be uniquely expressed as
%
%
\begin{eqnarray}
\label{eq211}
f=\sum_{k=0}^{\infty}\frac{H_k}{G^k}\;\;\; 
\text{ in  } K_{\infty}
\end{eqnarray}
where $H_0 \in \mathbb{A}$ and $H_k \in S_G\; (k \ge 1)$. 
Furthermore, $v_{\infty}(f)>0$ if and only if $H_0=0$. 
We call the expression (\ref{eq211}) the digit expansion of 
$f$ in base $G$. 
\label{prop21}
\end{prop}
%
%
\begin{proof}
Suppose $f=F_1/F_2 \; (F_1, F_2 \in \mathbb{A})$. 
Let $H_0, R \in \mathbb{A}$ with
 $F_1=H_0F_2+R$ and $\deg R<\deg F_2$. Putting $f_1=R/F_2$, then
 we have 
\[
f=H_0+f_1 \;\; \text { and } \;\;  v_{\infty}(f_1)>0. 
\]
For $k\ge 1$, we inductively  determine $H_k \in \mathbb{A}$
 and $f_{k+1} \in K$ satisfying with
\[
Gf_k=H_k+f_{k+1} \;\; \text { and } \;\;  v_{\infty}(f_{k+1})>0. 
\]
Then $H_k \in S_G\; (k\ge 1)$, and 
\begin{eqnarray*}
f=H_0+\frac{H_1}{G}+\cdots +\frac{H_k}{G^k}+\frac{f_{k+1}}{G^k}. 
\end{eqnarray*}
By $\lim_{k\rightarrow \infty}v_{\infty}(f_{k+1}/G^k)=\infty$, we 
get (\ref{eq211}). 

We next prove the uniqueness of (\ref{eq211}). Assume that
 \[
 \sum_{k=0}^{\infty}\frac{H_k}{G^k}=
 \sum_{k=0}^{\infty}\frac{I_k}{G^k},
 \]  
where $H_0, I_0 \in \mathbb{A}$ and
$H_k, I_k \in S_G\; (k\ge 1)$. 
Suppose that there exists an integer $k_0 \ge 0$ such
that $H_k=I_k\; (k<k_0)$ and $H_{k_0} \neq I_{k_0}$. Then
\begin{eqnarray*}
H_{k_0}-I_{k_0}=
\sum_{k=1}^{\infty}\frac{I_{k+k_0}-H_{k+k_0}}{G^k}.
\end{eqnarray*}
Since $v_{\infty}((I_{k+k_0}-H_{k+k_0})/G^k)>0$ for $k\geq 1$,
 we have $v_{\infty}(H_{k_0}-I_{k_0})>0$. 
This  contradicts $H_{k_0} -I_{k_0} \in \mathbb{A}\; \backslash \;\{0\}$.
Thus $H_k=I_k $ for $k \geq 0$.  

Finally, we see that 
\[
v_{\infty} (f)>0 \iff v_{\infty}(H_0)>0\iff H_0=0 \hspace{3mm}
\text {in (\ref{eq211})}
\]
since $v_{\infty}(H_k/G^k)>0$ for $k \geq 1$.
\end{proof}
For $M \in \mathbb{A}$, 
let $\mathcal{R}_M= \mathbb{A}/M  \mathbb{A}$. 
%
%
%
%
\begin{prop}
\label{prop22}
Let $G, M \in \mathbb{A}$ with $\deg G\geq 1, \; \deg M \geq 1$, 
 and $\gcd(G,M)=1$. 
Let $g$ be the order of $G \mod M$ in $ \mathcal{R}_M^{\times}$. 
Suppose that 
 $1/M=\sum_{k=1}^{\infty}H_k/G^k$ is
  the digit expansion of $1/M$
in base $G$. 
For  $k \geq 0$, we choose $G_k \in \mathbb{A}$ with
$\deg G_k <\deg M$ and $G_k \equiv G^k \pmod M$. 
\begin{itemize}
\item[(1)] For $k\geq 1$, we have
%
%
\begin{eqnarray}
\label{eq212}
H_k=\frac{GG_{k-1}-G_k}{M}. 
\end{eqnarray}
In particular, we obtain $H_k\neq 0 \; (k=1,2...)$ when 
$\deg G \geq \deg M$. 
\item[(2)] The sequence $\{H_k\}_{k\geq 1}$ is purely periodic, and
its period $l$ is equal to $g$.
\end{itemize}
\end{prop}
\begin{proof}
(1) Put $I_k=(GG_{k-1}-G_k)/M$ for $k\geq 1$. 
Since 
$
\deg (GG_{k-1}-G_{k})<\deg G+\deg M$ and $
GG_{k-1}\equiv G_k \pmod M$,
we have  $I_k \in S_G$. 
We see that 
\begin{eqnarray*}
\sum_{k=1}^{\infty}\frac{I_k}{G^k}
=\frac{1}{M}
\left\{\sum_{k=1}^{\infty}\frac{G_{k-1}}{G^{k-1}}- 
\sum_{k=1}^{\infty}\frac{G_{k}}{G^{k}}\right\}
=\frac{1}{M}=\sum_{k=1}^{\infty}\frac{H_k}{G^k}. 
\end{eqnarray*}
By the uniqueness of the digit expansion, 
we have $I_k=H_k$ for $k\geq 1$. 

(2) By assertion (1),  the sequence $\{H_k\}_{k\geq 1}$ is purely periodic, and 
$l \leq g$. On the other hand, we have
\begin{eqnarray*}
\frac{1}{M}=\sum^{l}_{k=1}\frac{H_k}{G^k}+\frac{1}{G^l}\cdot \frac{1}{M}. 
\end{eqnarray*}
Multiply both sides by $G^lM$ and consider modulo $M$. 
Then   $G^{l} \equiv 1 \pmod M$. We thus get $g \leq l$. 
\end{proof}
\begin{rem}
Equality (\ref{eq212}) plays an important role in 
the proof of our main results. For 
the number field version
of this equality, see \cite{gi1}, \cite{mi}, \cite{mt}. 
\end{rem}
%
%
\begin{cor}
\label{cor21}
Assume the notation in Proposition \ref{prop22}. 
Let $\alpha \in \mathbb{F}_q^{\times}$ of order $t$.
If $\gcd(M,G(\alpha G-1))=1$ and $t \mid g$, then 
\[
 \sum_{k=1}^{g}\alpha^kH_k=0. 
\]
\end{cor}
\begin{proof}
Put
\[
U=\sum_{k=1}^{g}\alpha^{k-1}G_{k-1}
 \; \left(=\sum_{k=1}^{g}\alpha^{k}G_{k}\right).
\]
By $\gcd(M,G(\alpha G-1))=1$ and $t \mid g$, we have 
\[
U \equiv \sum_{k=1}^{g}\alpha^{k-1}G^{k-1}=\frac{(\alpha G)^g-1}{\alpha G-1} 
\equiv 0 \pmod M . 
\]
This gives $U=0$
since $\deg U<\deg M$. 
Therefore, by (\ref{eq212}), we obtain
\[
\sum_{k=1}^{g}\alpha^kH_k=\frac{(\alpha G-1)U}{M}=0. 
\]
\end{proof}
%
%
\begin{rem}
Considering the case $\alpha=1$ in Corollary \ref{cor21}, 
we obtain
\[
 \sum_{k=1}^{g}H_k=0. 
\]
This result was discovered by Rudnick \cite{ru}. 
\end{rem}

%
%
%
%
\subsection{Cyclotomic function fields and divisor class numbers}
\quad 
In this subsection, we review  basic facts about cyclotomic
function fields and divisor class numbers. 
For details, see \cite{gr, ha, ro}. 

Let $K^{ac}$ be an algebraic closure of $K$. 
We define the action
\[
I*x=I(\varphi+\mu)(x) \;\;\;\; (x  \in K^{ac}, \; I \in \mathbb{A}),
\]
where $\varphi$ and  $\mu$ are the $\mathbb{F}_q$-linear maps 
defined by
\begin{eqnarray*}
\varphi:  K^{ac} \rightarrow K^{ac}\;\;\; x \mapsto x^q, \hspace{8mm}
\mu:  K^{ac} \rightarrow K^{ac}\;\;\; x \mapsto Tx. 
\end{eqnarray*}
By the above action, $K^{ac}$ becomes an 
$\mathbb{A}$-module.
Suppose that $P \in \mathbb{A}$ is 
 monic irreducible of degree $d$. 
Put
\[
\Lambda_P=\{x \in K^{ac}\; |\; P*x=0\}.
\]
This is a cyclic $\mathbb{A}$-submodule of $K^{ac}$. 
Fix a generator $\lambda_P$ of $\Lambda_P$. 
Then we have the following isomorphism of $\mathbb{A}$-modules
\[
\mathcal{R}_P\longrightarrow \Lambda_P
\;\;\; (\overline{I} \rightarrow I*\lambda_P),
\]
where 
$\mathcal{R}_P=\mathbb{A}/P\mathbb{A}$, and 
$\overline{I} =I \text{ mod } P \in \mathcal{R}_P$. 
Let $K_P=K(\Lambda_P)$. Then $K_P/K$ is a Galois extension, 
and we have  the following isomorphism
\begin{eqnarray*}
\mathcal{R}_P^{\times} \rightarrow \text{Gal} (K_P/K)
\;\;\; (\overline{I} \mapsto \sigma_{I}),
\end{eqnarray*}
where 
$\text{Gal} (K_P/K)$ is the Galois group of $K_P/K$, and 
$\sigma_{I }$  is the isomorphism given
by $\sigma_{I }(\lambda_P)=I*\lambda_P$. 
Since 
$\mathcal{R}_P$ is a finite field of order $q^d$, 
the extension $K_P/K$ is 
 cyclic of degree $q^d-1$.
The field $K_P$ is called the $P$-th cyclotomic function field. 


Let $L$ be an intermediate field of $K_P/K$. 
We denote by $H_L$ the
subgroup of $\mathcal{R}_P^{\times}$
corresponding to $L$. 
The field $L^+=L\cap K_\infty$ is called the maximal real subfield of $L$.
By Theorem 12.14 in \cite{ro}, we have 
$H_{L^+}=H_L\mathbb{F}_q^{\times}$.

The next lemma is well-known to the experts.
But we give a proof for the reader's convenience. 
%
%
%
%
\begin{lem}
\label{lem21}
Suppose that $q$ is odd, and 
$P \in \mathbb{A}$ is monic irreducible of degree $d$. 
Then  $L=K\left(\sqrt{(-1)^dP}\right)$
is contained in $ K_P$. Furthermore, 
%
%
\begin{eqnarray}
\label{eq221}
 L^+=
 \begin{cases}
L
 \hspace{11mm} \text{if $d$ is even, }  \\[1mm]
 K
   \hspace{10mm}\text{if $d$ is odd}. \end{cases}
\end{eqnarray}
\end{lem}
\begin{proof}
By the discussion of the proof of Lemma 16.13 in \cite{ro},
we can take an element $\pi \in K_P$ satisfying with 
$
\pi^{q-1} =(-1)^rP,
$
where $r=(q^d-1)/(q-1)$. 
Hence $L$ is contained in $K_P$
 since $r \equiv d  \pmod 2$ and $q$ is odd.

We note that $K_P/K$ is cyclic and 
$
 [K_P^+:K]=r\equiv d \pmod 2
 $. Therefore,  
 \[
L=L^+ \iff  L \subseteq K_P^+ \iff d\equiv 0 \pmod 2. 
 \]
This implies (\ref{eq221}). 
\end{proof}
Let $X_P$ be the group of Dirichlet characters modulo $P$. 
For an intermediate field $L$ of $K_P/K$, we set
\begin{eqnarray*}
X_L &=&\{\chi \in X_P \; |\; \chi\left(\overline{I}\right)=1 
\;\; (\overline{I} \in H_L)\}, \\[2mm]
X_L^+&=&\{\chi \in X_P \; |\; \chi\left(\overline{I}\right)=1 \;\;
(\overline{I} \in H_{L^+})\}, 
\end{eqnarray*}
and $X_L^-=X_L \; \backslash \;X_L^+$. 
We denote by $h_L$ (resp. $h_L^+$ )  the divisor class number of $L$
(resp. $L^+$). 
By a similar argument as in the proof of
 Theorem 2 in \cite{gr}, we have 
\begin{eqnarray}
&& 
\label{eq222}
h_L=\prod_{\chi \in X_L^+\atop \chi \neq 1_P}  
\left\{-\sum_{I \in \mathbb{M} \atop \deg I<d}\chi(I) \deg I\right\}
\; \cdot \; \prod_{\chi \in X_L^-}  
\left\{\sum_{I \in \mathbb{M} \atop \deg I<d}\chi(I) \right\},
\\[2mm]
 &&
 \label{eq223}
  h_L^+=
  \prod_{\chi \in X_L^+\atop \chi \neq 1_P}  
\left\{-\sum_{I \in \mathbb{M} \atop \deg I<d}\chi(I) \deg I\right\}. 
\end{eqnarray}
Here, 
$\mathbb{M}$ is the set of all monic
polynomials of $\mathbb{A}$,  
and $1_P$ is the 
trivial character of $X_P$. 
The rational number $h_L^-=h_L/h_L^+$ is called
the relative divisor class number of $L$.
 By (\ref{eq222}) and (\ref{eq223}), we have 
\begin{eqnarray}
\label{eq224}
h_L^-=\prod_{\chi \in X_L^-}  
\left\{\sum_{I \in \mathbb{M} \atop \deg I<d}\chi(I) \right\}. 
\end{eqnarray}
The right-hand side of (\ref{eq224}) is an algebraic integer. Hence
the rational number $h_L^-$ is an integer. 

%
%
%
%
%
%
%
%
\section{Main results }
Suppose that $P \in \mathbb{A}$ is  monic irreducible \vspace{1mm}
and $G \in \mathbb{A}$ is a primitive root modulo $P$. 
Let
$
d=\deg P, \;  e=\deg G, \;  r=(q^d-1)/(q-1). 
$
Let
$
1/P=\sum^{\infty}_{k=1}H_k/G^k
$\vspace{1mm}
be the digit expansion of $1/P$ in base $G$. 
\vspace{1mm}
Choose $G_k \in \mathbb{A} \; (0 \leq k \leq q^d-1)$ with
\vspace{1mm}
$\deg G_k<d$ and $G_k \equiv G^k \pmod P$. 
%
%
%
 %
\begin{lem}
\label{lem31}
For $ 0\leq s \leq d-1$, we have
\begin{eqnarray*}
\left\{ \left[\hspace{1pt}\overline{I} \hspace{1pt} \right]\; \Big|\; 
\begin{array}{ll}I \in \mathbb{M},\\[1mm]
 \deg I=s
 \end{array}\right\}
=\left\{ \left[\hspace{1pt}\overline{G_{k}} \hspace{1pt}\right]\; \Big|\;
\begin{array}{ll}
 0\leq k\leq r-1, \\[1mm]
  \deg G_{k}=s
 \end{array} \right\},
\end{eqnarray*}
where $[\hspace{1pt}\overline{I} \hspace{1pt}]$ is the 
equivalent class of $\overline{I} \in \mathcal{R}_P$ in 
$\mathcal{R}_P^{\times}/\mathbb{F}_q^{\times}$. 
\end{lem}
\begin{proof}
Since $G$ is a primitive root of modulo $P$, we have 
\[
\overline{G_k}=(\overline{G})^k \in \mathbb{F}_q^{\times}
\iff r\; | \;  k. 
\]
Hence 
$\{\; \overline{G_k}\; |\; 0 \leq k \leq r-1\; \}$ is a complete residue
system of $\mathcal{R}_P^{\times}/\mathbb{F}_q^{\times}$. 
This implies Lemma \ref{lem31}. 
\end{proof}

We first consider
the plus-part of divisor class numbers. 
We put
\begin{eqnarray*}
F^{(+)}(u)=\sum_{k=1}^{r}(\deg H_k )u^{k-1}.
\end{eqnarray*}
In the above definition,  
we understand $\deg H_k$ as $0$ if $H_k=0$. 
Let
\[
X_P^+=\{\chi \in X_P\; |\; \chi(a)=1\;\text{ for all }
 a \in \mathbb{F}_q^{\times} \}.
\]
Then we have
%
%
%
%
\begin{prop}
\label{prop31}
For $\chi \in X_P^+$, we have 
\[
F^{(+)}(\zeta)=\sum_{I \in \mathbb{M} \atop d-e \leq \deg I<d}
(\deg I -(d-e))\chi(I), \\[2mm]
\]
where $\zeta=\chi(G)$. 
\end{prop}
%
%
\begin{proof}
By (\ref{eq212}), we have $GG_{k-1}=H_kP+G_k$. It follows that 
\begin{eqnarray*}
\deg H_k=\begin{cases}
\deg G_{k-1}-(d-e) \hspace{9mm} \text{if $\deg G_{k-1}\geq d-e$, }\\[2mm]
0 \hspace{36mm} \text{if $\deg G_{k-1}<d-e$.}
\end{cases}
\end{eqnarray*}
Define
$
\text{Deg}:\mathcal{R}_P^{\times}\rightarrow \mathbb{Z} \; 
(\overline{I}\mapsto \deg r_I), 
$
where $r_I$ is the remainder of $I$
divided by $P$. 
It is easy to check that  $\text{Deg}$  and $\chi$ are both functions on 
$\mathcal{R}_P^{\times}/\mathbb{F}_q^{\times}$. 
Therefore, by Lemma \ref{lem31},  we have %
\begin{eqnarray*}
&&\sum_{I \in \mathbb{M} \atop d-e \leq \deg I<d} 
(\deg I -(d-e))\chi(I) \\
&& \hspace{15mm}=
\sum_{I \in \mathbb{M} \atop d-e \leq \deg I<d}
(\text{Deg}([\hspace{1pt}\overline{I}\hspace{1pt}])
 -(d-e))\chi([\hspace{1pt}\overline{I}\hspace{1pt}]) \\[2mm]
&& \hspace{15mm} =
\sum_{1 \leq k\leq r \atop d-e \leq \deg G_{k-1}}
(\text{Deg}([\hspace{1pt} \overline{G_{k-1}}\hspace{1pt}]) -(d-e))
\chi( [\hspace{1pt} \overline{G_{k-1}}\hspace{1pt}]) \\[2mm]
&&  \hspace{15mm}= \sum_{k=1}^r (\deg H_k) \zeta^{k-1}. 
\end{eqnarray*}
We thus get Proposition \ref{prop31}. 
\end{proof}
%
%
%
%
\begin{cor}
\begin{eqnarray*}
\sum_{k=1}^{r}\deg H_k=
\begin{cases}
q^{d-e}s_1(e) \hspace{24mm}\text{if $e<d$}, \\[2mm]
(e-d)s_0(d)+s_1(d) \hspace{7.5mm} \text{if $e\geq d$,}
\end{cases}
\end{eqnarray*}
where $s_0(j)=\sum_{i=0}^{j-1}q^i$ and $s_1(j)=\sum_{i=0}^{j-1}iq^i$.
\end{cor}
\begin{proof}
We consider the case $\chi=1_P$
 in Proposition  \ref{prop31}.  
 Assume that  $e<d$. 
Since the number of monic polynomials 
of degree $i$ in $\mathbb{A}$ is equal to $q^i$, we have 
\begin{eqnarray*}
\sum_{k=1}^{r}\deg H_k
&=& \sum_{i=0}^{e-1}\sum_{I \in \mathbb{M} \atop  \deg I=d-e+i}
(\deg I-(d-e)) \\
&=&\sum_{i=0}^{e-1}i q^{d-e+i} \\[1mm]
&=& q^{d-e}s_1(e). 
\end{eqnarray*}
We next assume that  $e\ge  d$. Then
\begin{eqnarray*}
\sum_{k=1}^{r}\deg H_k 
&=& \sum_{i=0}^{d-1}\sum_{I \in \mathbb{M} \atop  \deg I=i}
(\deg I-(d-e)) \\
&=&\sum_{i=0}^{d-1} (i+(e-d))q^{i} \\[1mm]
&=& s_1(d)+(e-d)s_0(d). 
\end{eqnarray*}
\end{proof}
%
%
%
%
\begin{cor}
\label{cor32}
If $\chi \in X_P^+ \; \backslash \;\{1_P\}$ and $e\geq d$, 
then 
\begin{eqnarray*}
F^{(+)}(\zeta)=\sum_{I \in \mathbb{M} \atop \deg I<d} \chi(I)\deg I,
\end{eqnarray*}
where $\zeta=\chi(G)$. 
\end{cor}
%
%
\begin{proof}
Since $\chi(\mathbb{F}_q^{\times})=\{1\}$, we have
\begin{eqnarray*}
\sum_{I \in \mathbb{M} \atop \deg I<d} \chi(I)
=\frac{1}{q-1}\sum_{I \in \mathbb{A} \atop \deg I<d} \chi(I)
=0. 
\end{eqnarray*}
Therefore, by Proposition \ref{prop31}, we obtain
\begin{eqnarray*}
F^{(+)}(\zeta)=\sum_{I \in \mathbb{M} \atop  \deg I<d}
\chi(I)(\deg I -(d-e))=
\sum_{I \in \mathbb{M} \atop \deg I<d}
\chi(I)\deg I. 
\end{eqnarray*}
\end{proof}

%
%
%
%
\begin{thm}
\label{th31}
Let $L$ be an intermediate field of $K_P/K$. 
Assume that  $m=[L^+:K]>1$ and $e \geq d$. Then,
\begin{eqnarray*}
(-1)^{m-1}h_L^+=
\prod_{\zeta \in U_L^+\atop \zeta \neq 1}F^{(+)}(\zeta),
\end{eqnarray*}
where $U_L^+$ is the set of all $m$-th roots of unity.  
\end{thm}
\begin{proof}
Since $X_L^+$ is 
 a cyclic group of order $m$, the map
$
X_L^+\longrightarrow U_L^+ \; (\chi \mapsto \chi(G))
$ is a group isomorphism. 
By (6) and Corollary \ref{cor32},  we have
\begin{eqnarray*}
(-1)^{m-1}h_L^+
&=&\prod_{\chi \in X_L^+\atop \chi \neq 1_P} 
\left\{\sum_{I \in \mathbb{M} \atop \deg I<d}\chi(I)\deg I \right\}\\[2mm]
&=&\prod_{\zeta \in U_L^+ \atop \zeta \neq 1}F^{(+)}(\zeta).
\end{eqnarray*}
\end{proof}
Now we prove the even case of Theorem \ref{th12}.
%
%
%
%
\begin{proof}
By Lemma \ref{lem21}, the field
$L$ is contained in $K_P$ and $L=L^+$. 
Therefore, by Theorem \ref{th31}, we obtain
\begin{eqnarray*}
h_L=h_L^+=-F^{(+)}(-1)=
\sum_{k=1}^{r}(-1)^k\deg H_k. 
\end{eqnarray*}
\end{proof}
%
%
%
%
%
%
\begin{ex}
\label{ex31}
Suppose that $q=3$, $P=T^2+1 $, $G=T^2+T+2$, and 
$L=K(\sqrt{P})$. 
Then $P$ is irreducible in $\mathbb{F}_3[T]$, and
$G$ is a primitive root modulo $P$. Let
$1/P=\sum_{k=1}^{\infty}H_k/G^k$ be
the digit expansion of $1/P$ in base $G$.
Then $r=(q^2-1)/(q-1)=4$, and 
\[
H_1=1, \;\;\; 
H_2=T+2, \;\;\; 
H_3=2T+2, \;\;\; 
H_4=2T. 
\]
By Theorem \ref{th12}, we have
\[
h_{L}=(\deg H_2+\deg H_4)-
(\deg H_1+\deg H_3)=1.
\]
We next calculate 
the divisor class number of $K_P^+$.
 We see that $[K_P^+:K]=4$, and   
\begin{eqnarray*}
F^{(+)}(u) =\sum_{k=1}^{4}(\deg H_k)u^{k-1}=u+u^2+u^3, \;\;\; 
U_{K_P}^+=\{\zeta \in \mathbb{C}\; |\; \zeta^4=1\}. 
\end{eqnarray*}
By Theorem \ref{th31},  we have
\[
h_{K_P}^+=(-1)^{[K_P^+:K]-1}
\prod_{\zeta \in U_{K_P}^+ \atop \zeta \neq 1}F^+(\zeta)=1. 
\]
\end{ex}
%
%
\begin{ex}
Suppose that $q=2$, $P=T^3+T+1$, $G=T^3$. 
Then $P$ is irreducible in $\mathbb{F}_2[T]$, and
$G$ is a primitive root modulo $P$. Let
$1/P=\sum_{k=1}^{\infty}H_k/G^k$ be
the digit expansion of $1/P$ in base $G$. Then
$r=(q^3-1)/(q-1)=7$, and 
\begin{eqnarray*}
&&H_1=1, \;\;\; 
H_2=T+1, \;\;\; 
H_3=T^2, \;\;\; 
H_4=T^2+1, \\
&&
H_5=T^2+T, \;\;\; 
H_6=T, \;\;\; 
H_7=T^2+T+1. 
\end{eqnarray*}
We see that $[K_P^+:K]=7$, 
and 
\begin{eqnarray*}
F^{(+)}(u) &=& \sum_{k=1}^{7}(\deg H_k)u^{k-1}=
u+2u^2+2u^3+2u^4+u^5+2u^6, \\[2mm]
U_{K_P}^+&=& \{\zeta \in \mathbb{C}\; |\; \zeta^{7}=1\}. 
\end{eqnarray*}
By Theorem \ref{th31},  we have
\[
h_{K_P}^+=(-1)^{[K_P^+:K]-1}
\prod_{\zeta \in U_{K_P}^+\atop \zeta \neq 1}F^+(\zeta)=71. 
\]
\end{ex}

%
%
%
%
%
We next consider the minus-part
of divisor class numbers.
Denote by $\widehat{\mathbb{F}_q^{\times}}$
the character group of $\mathbb{F}_q^{\times}$. 
Recall that $\delta(I)$ is the leading coefficient of 
 $I \in \mathbb{A}$ if $I \neq 0$, 
 and $\delta(0)=0$.  
 Suppose that $\lambda \in \widehat{\mathbb{F}_q^{\times}}$. 
We define
\begin{eqnarray*}
F^{(\lambda)}(u)=
\lambda(\delta(G))\sum_{k=1}^{r}\bar{\lambda}(\delta(H_k))u^{k-1}, 
\end{eqnarray*}
where $\bar{\lambda} $ is
 the complex conjugate of  $\lambda$. 
In the above definition, we understand $\bar{\lambda}(\delta(H_k))$
as $0$ if $\delta(H_k)=0$. 
We set
\begin{eqnarray*}
X_P^{(\lambda)} 
= \{ \chi \in X_P\; |\; \chi|_{\mathbb{F}_q^{\times}}=\lambda \},
\end{eqnarray*}
where $\chi|_{\mathbb{F}_q^{\times}}$ is the restriction to 
$\mathbb{F}_q^{\times}$ of $\chi$.  
%
%
%
%
\begin{prop}
\label{prop32}
For $\chi \in X_P^{(\lambda)}$, we have
\[
F^{(\lambda)}(\zeta)=\sum_{I \in \mathbb{M} \atop d-e \leq \deg I<d}
\chi(I),
\]
where $\zeta=\chi(G)$. 
\end{prop}
%
%
\begin{proof}
Define $H_{\chi} :
\mathcal{R}_P^{\times}\rightarrow \mathbb{C}$ by
\begin{eqnarray*}
H_{\chi}(\overline{I})=
\overline{\lambda}(\delta(r_I))\chi(r_I), 
\end{eqnarray*}
where $r_I$ is the remainder of $I$ divided by $P$. 
It is easy to check that $H_{\lambda}$ is a function on 
$\mathcal{R}_P^{\times}/\mathbb{F}_q^{\times}$. 
By Lemma \ref{lem31},  we have 
\begin{eqnarray*}
\sum_{I \in \mathbb{M} \atop d-e \leq \deg I<d}\chi(I)
&=& \sum_{I \in \mathbb{M} \atop d-e \leq \deg I<d}
\overline{\lambda}(\delta(r_I))\chi(r_I) \\[2mm]
&=& \sum_{I \in \mathbb{M} \atop d-e \leq \deg I<d}
H_{\chi}([\hspace{1pt}\overline{I}\hspace{1pt}] ) \\[2mm]
&=& \sum_{1 \leq k\leq r \atop d-e \leq \deg G_{k-1} }
H_{\chi}([\hspace{1pt}\overline{G_{k-1}} \hspace{1pt}]) \\[2mm]
&=&  \sum_{1 \leq k\leq r \atop d-e \leq \deg G_{k-1}}
\overline{\lambda}(\delta(G_{k-1}))\zeta^{k-1}.
\end{eqnarray*}
By $GG_{k-1}=H_kP+G_k$, we have
\begin{eqnarray*}
\delta(H_k)=\begin{cases}
\delta(G)\delta(G_{k-1}) \hspace{10.5mm} \text{if $\deg G_{k-1}\geq d-e$, }\\[2mm]
0 \hspace{28.5mm} \text{if $\deg G_{k-1}<d-e$.}
\end{cases}
\end{eqnarray*}
Therefore, 
\begin{eqnarray*}
\sum_{I \in \mathbb{M} \atop d-e \leq \deg I<d}\chi(I)
=\lambda(\delta(G))
 \sum_{k=1}^r
\overline{\lambda}(\delta(H_k))\zeta^{k-1} =F^{(\lambda)}(\zeta). 
\end{eqnarray*}
\end{proof}

%
%
%
%
\begin{cor}
\label{cor33}
Assume that $\chi \in X_P^{(\lambda)}$ and $e\geq d$. Then, 
\[
F^{(\lambda)}(\zeta)=\sum_{I \in \mathbb{M}\atop \deg I<d}\chi(I).
\]
\end{cor}
\begin{proof}
This follows immediately from Proposition \ref{prop32}. 
 \end{proof}

Let $L$ be an intermediate field of $K_P/K$. 
Put
\[
l=[L:K],\;\;  m=[L^+:K], \;\; n=[L:L^+].
\] 
Since $X_P$ is cyclic, we have
\begin{eqnarray*}
X_L =\{\chi \in X_P \; | \; \chi^l=1_P\}, \;\;\;\; 
X_L^+ =\{\chi \in X_P \; | \; \chi^m=1_P\}.
\end{eqnarray*}
We set
\begin{eqnarray*}
Y_L=\{\lambda \in  \widehat{\mathbb{F}_q^{\times}}
\;| \; \lambda^n=\lambda_0
\}, \;\;\; U_n=\{\zeta \in \mathbb{C} \; |\; \zeta^n=1\}, 
\end{eqnarray*}
where $\lambda_0$ is the trivial character of $\mathbb{F}_q^{\times}$. 
Then the maps
\begin{eqnarray*}
&& \varphi : X_L/X_L^+ \longrightarrow Y_L\; 
(\chi X_L^+\mapsto \chi |_{\mathbb{F}_q^{\times}}), \; \;\; \\[2mm]
&& \psi : X_L/X_L^+ \longrightarrow U_n\; 
(\chi X_L^+\mapsto \chi(G)^m)
\end{eqnarray*}
are both group isomorphisms. 
For $\lambda \in Y_L$, we put
 $\alpha_{\lambda}=\psi\circ \varphi^{-1}(\lambda)$. 
%
%
%
%
\begin{thm}
\label{th32}
Let $L$ be an intermediate field of $K_P/K$. Assume that
$L \neq L^+$ and $e\ge d$. Then 
\begin{eqnarray*}
h_L^-=\prod_{\lambda \in Y_L \atop \lambda \neq \lambda_0}
\prod_{\zeta \in U_L^{(\lambda)}}F^{(\lambda)}(\zeta),
\end{eqnarray*}
where 
$ U_L^{(\lambda)}=\{\zeta \in \mathbb{C}\; |\; \zeta^m=\alpha_\lambda
\}$.
\end{thm}
\begin{proof}
For $\lambda \in Y_L$, we set
 \begin{eqnarray*}
X_L^{(\lambda)} 
= \{\chi \in X_L\; |\; \chi|_{\mathbb{F}_q^{\times}}=\lambda\}.
\end{eqnarray*}
Then,
\begin{eqnarray*}
X_L^-=X_L \;\backslash \; X_L^+=
\bigcup_{\lambda \in Y_L \atop \lambda \neq \lambda_0} X_L^{(\lambda)}.
\end{eqnarray*} 
By the definition of $\alpha_{\lambda}$, 
we have
\[
X_L^{(\lambda)}= \{\chi \in X_L\; |\; \chi(G)^m=\alpha_\lambda\}. 
\]
This implies $U_L^{(\lambda)}=\{\chi(G)\; |\; \chi \in X_L^{(\lambda)}\}$. 
Therefore, by (\ref{eq224}) and Corollary \ref{cor33},  we have 
\begin{eqnarray*}
h_L^-
&=&\prod_{\chi \in X_L^-} 
\left\{\sum_{I \in \mathbb{M} \atop \deg I<d}
\chi(I)\right\}\\[2mm]
&=&
\prod_{\lambda \in Y_L\atop \lambda \neq \lambda_0}
\prod_{\chi \in X_L^{(\lambda)}}
\left\{\sum_{I \in \mathbb{M}\atop \deg I<d}\chi(I)\right\} \\[2mm ]
&=&
\prod_{\lambda \in Y_L\atop \lambda \neq \lambda_0}
\prod_{\zeta \in U_L^{(\lambda)}}F^{(\lambda)}(\zeta).
\end{eqnarray*}
\end{proof}

%
%
%
We prove the odd case of Theorem \ref{th12}. 
%
%
\begin{proof}
By Lemma \ref{lem21}, the field $L$
is contained in $K_P$ and  $L^+=K$. 
Suppose that $\lambda$ is 
the quadratic character of $\mathbb{F}_q^{\times}$. 
Then  $Y_L$ is generated by 
$\lambda$. 
Note that  $\alpha_{\lambda}=-1$ and 
\begin{eqnarray*}
U_L^{(\lambda)}=\{-1\}, \hspace{5mm}
\lambda(\delta(G))=\eta_G, \hspace{5mm}
\overline{\lambda}(\delta(H_k))
=\lambda(\delta(H_k))=\varepsilon_k.
\end{eqnarray*}
Therefore, by Theorem \ref{th32}, we obtain
\begin{eqnarray*}
h_L=h_L^-=
F^{(\lambda)}(-1)=
(-\eta_G)\sum_{k=1}^{r}(-1)^k\varepsilon_k. 
\end{eqnarray*}
\end{proof}
%
%
\begin{ex}
Suppose that $q=3$, $P=T^3+2T+2$, $G=T^3+T+2$, and 
$L=K(\sqrt{-P})$. 
Then $P$ is irreducible in $\mathbb{F}_3[T]$, and
$G$ is a primitive root modulo $P$. 
Let 
$1/P=\sum_{k=1}^{\infty}H_k/G^k$ be
the digit expansion of $1/P$ in base $G$. Then
$r=(q^3-1)/(q-1)=13$, $\eta_G=1$, and 
\begin{eqnarray*}
&&H_1=1, \;\;\; 
H_2=2T, \;\;\; 
H_3=T^2+2, \;\;\; 
H_4=2T+2,\;\;\; 
H_5=T^2+T+2, \\
&& 
H_6=2T^2+2T, \;\;\; 
H_7=T^2+2T, \;\;\; 
H_8=T^2+T+1, \;\;\; 
H_9=2T^2, \\
&& 
H_{10}=2T+1, \;\;\; 
H_{11}=T^2+2T+2, \;\;\; 
H_{12}=T^2+2T+1,\;\;\;
H_{13}=T^2+1.  
\end{eqnarray*}
Let $\varepsilon_k=1 \text{ or } -1 \; (k=1,2, ..., 13)$
 according to $\delta(H_k)=1$ or $\delta(H_k)=2$. 
By Theorem \ref{th12}, we have
\[
h_L=(-\eta_G)\sum_{k=1}^{13}(-1)^k\varepsilon_k=7. 
\]
We next calculate the relative class number of $K_P$.
Let $\lambda$
be the quadratic character of $\mathbb{F}_3^{\times}$. 
Then $[K_P^+:K]=13$, $\alpha_{\lambda}=-1$, and 
\begin{eqnarray*}
F^{(\lambda)}(u) &=&\lambda(\delta(G))\sum_{k=1}^{13}\overline{\lambda}
(\delta(H_k))u^{k-1}  \\
&=&
1-u+u^2-u^3+u^4-u^5+u^6+u^7-u^8-u^9+u^{10}+u^{11}+u^{12}, \\[2mm]
U_{K_P}^{(\lambda)}&=&\{\zeta \in \mathbb{C}\; |\; \zeta^{13}=-1\}. 
\end{eqnarray*}
Therefore, by Theorem \ref{th32},  we have
\[
h_{K_P}^-=\prod_{\zeta \in U_{K_P}^{(\lambda)}}F^{(\lambda)}(\zeta)
=774144=2^{12}\cdot 3^{3}\cdot 7.
\]
\end{ex}
%
%
%
%

%
%
%
\text{ }\\
Daisuke Shiomi \\
Department of Science \\
Faculty of Science, Yamagata University \\
Kojirakawa-machi 1-4-12 \\
Yamagata 990-8560, Japan \\
shiomi@sci.kj.yamagata-u.ac.jp

\end{document}